\documentclass[preprint,12pt]{elsarticle}

\usepackage{amssymb}

\usepackage{graphicx}
\usepackage{amssymb}
\usepackage{newlfont}
\usepackage{t1enc}
\usepackage{amsfonts}

\newtheorem{thm}{Theorem}[section]

\newproof{pf}{Proof}

\newcommand{\A}{\mathcal{A}}
\let \dis =\displaystyle

\usepackage{relsize}
\usepackage{amssymb}

\begin{document}

\begin{frontmatter}



\title{THE REAL POWERS OF THE CONVOLUTION
OF A GAMMA DISTRIBUTION AND A BERNOULLI DISTRIBUTION}

\author{Ben Salah Nahla }
\ead{bensalahnahla@yahoo.f}
\author{Masmoudi Afif }
\ead{Afif.Masmoudi@fss.rnu.tn}

\cortext[cor1]{Corresponding author} \cortext[cor2]{Principal
corresponding author}
\address{Laboratory of Probability and  Statistics, Faculty
of Science of Sfax, Sfax University. B. P. 1171\\ Sfax. Tunisia}

\begin{abstract}
In this paper, we essentially compute  the set of $x,y>0$ such
that the mapping $z \longmapsto \Big{(}1-r+r e^z\Big{)}^x
\Big{(}\dis\frac{\lambda}{\lambda-z}\Big{)}^{y}$ is a Laplace
transform. If $X$ and $Y$ are two independent random variables
which have respectively Bernoulli and Gamma distributions, we
denote by $\mu$ the distribution of $X+Y.$ The above problem is
equivalent to finding the set of $x>0$ such that $\mu^{{\ast}x}$
exists.

\end{abstract}


\end{frontmatter}


\section{Introduction and Preliminaries}\label{intro}

We introduce first some notations and review some basic concepts
concerning the Laplace transform and J{\o}rgensen set. For more
details, we refer the reader to \cite{B, JOR}. \\ For a positive
measure $\mu$ on $\mathbb{R}$, we denote by
$$L_{\mu}:\mathbb{R}\longrightarrow (0;+\infty); \theta \longmapsto \int_{\mathbb{R}} e^{\theta x} \mu(dx)$$
its Laplace transform, and let
$$\Theta(\mu)=\mathrm{interior}\{\theta; L_{\mu}(\theta)< +\infty \}.$$
We denote by $\mathcal{M}(\mathbb{R})$ the set of such measures on
$\mathbb{R}$ such that, furthermore, $\mu$ is not concentrated on
a point and such that $\Theta(\mu)$ is not empty. \newpage
\noindent For $\mu \in
\mathcal{M}(\mathbb{R}),$  $L_\mu$ is real analytic and strictly
convex on $\Theta(\mu)$.

The J{\o}rgensen parameter is the parameter corresponding to the
power of convolution such that it is the variance in the Gaussian
distribution and it is the shape parameter in the Gamma distribution.

Let $\mu$ be an element of $\mathcal{M}(\mathbb{R}),$ the
J{\o}rgensen set $\Lambda(\mu)$ of $\mu$ is the set of $x > 0$
such that there exists $\mu_x \in \mathcal{M}(\mathbb{R})$ with
$\Theta(\mu_x) = \Theta(\mu)$ and $L_{\mu_{x}}(\theta) =
(L_{\mu}(\theta))^{x}$ (see \cite{JOR}).

In this case, $\mu_{x}$ is called the $x^{th}$ convolution power of
$\mu.$

If $x$ and $x'$ are two elements of the J{\o}rgensen set
$\Lambda(\mu)$ of $\mu,$ then $(\mu_{x},\mu_{x'})$ is convolvable
since $\Theta(\mu_{x})=\Theta(\mu_{x'}).$ Furthermore $x+x' \in
\Lambda(\mu).$ Hence, $\Lambda(\mu)$ is a semigroup under addition
and
$$\mu_{x+x'}=\mu_{x}\ast \mu_{x'}.$$ It
contains 1 by definition, therefore it contains the set
$\mathbb{N}^{\ast}$ of non-negative integers. The calculation of
$\Lambda(\mu)$ is sometimes a hard problem: it is
$\mathbb{N}^{\ast}$ when $\mu$ is the  Bernoulli distribution on
$\{0, 1\}$.
\\A probability distribution $\mu$ on $\mathbb{R}$ is
infinitely divisible if, for every integer $n,$ there exists a
distribution $\mu_n$ such that
$$\mu=\mu_n^{\ast n}$$ that is, $\mu$ is the $n^{th}$ power of
convolution of $\mu_n.$

In other words, $\mu$ is infinitely divisible if, for each integer
$n,$ it can be represented as the distribution of the sum
$S_n=X_{1,n}+X_{2,n}+....+X_{n,n}$ of $n$ independent random
variables with a common distribution $\mu_n$ (see \cite{feller}).

If $\mu$ is an element of $\mathcal{M}(\mathbb{R}),$ then $\mu$ is
infinitely divisible if and only if its J{\o}rgensen set
$\Lambda(\mu)$ is equal to $(0,+\infty).$
\\Let $\nu$ be a distribution on the real line
having Laplace transform and which is not infinitely divisible.
Consider now an infinitely divisible distribution $\nu'$ on the
real line, also having Laplace transform. We denote by $\nu'_{y}$
the distribution such that $L_{\nu'_y}={L^{y}_{\nu'}}.$
\\Let $\mu=\nu \ast \nu'$ be the convolution product of $\nu$ and
$\nu'.$
\\Letac et al. \cite{R} considered the case where $\nu$ is
Bernoulli  and $\nu'$ is negative Binomial. In this case, the
problem is equivalent to finding the set of $(x,y) \in
(0,+\infty)^{2}$ such that there exists a probability $\mu_{x,y}$
on the real line with Laplace transform $L_{\nu}^{x}L_{\nu'}^{y}.$
\\In this work, we consider the case where $\nu$ is
Bernoulli and $\nu'$ is Gamma. In this situation, the techniques
are completely different since  we use essentially the analyticity
of the Laplace transform.
\\ For fixed $r \in (0,1)$ and $\lambda>0,$ the present paper
wonders for which values of $x,y>0$ the function defined on
$(-\infty,\lambda)$ by $z \longmapsto \Big{(}1-r+r e^z\Big{)}^x
\Big{(}\dis\frac{\lambda}{\lambda-z}\Big{)}^{y}$ is the Laplace
transform of a probability. This function is the Laplace transform
of the function $$t \longmapsto
\frac{\lambda^{y}(1-r)^{x}}{\Gamma(y)} \dis\sum_{k \geq
0}\frac{x(x-1)(x-2)....(x-k+1)}{k!}
\Big{(}\frac{r}{1-r}\Big{)}^{k} e^{-\lambda(t-k)}(t-k)_{+}^{y-1}$$
where $a_{+}$ means $\max(a,0)$. Denote $R=\frac{r}{1-r}
e^{\lambda}$ for simplicity. The problem is therefore equivalent
in finding the set $\Lambda_{R}$ of the $(x,y)'s$ such that
$$t \longmapsto
\frac{\lambda^{y}(1-r)^{x}}{\Gamma(y)} \dis\sum_{k \geq
0}\frac{x(x-1)(x-2)....(x-k+1)}{k!} R^{k} e^{-\lambda
t}(t-k)_{+}^{y-1}=f(t)$$ is a positive function for all $t$ (note
that this series converges, having a general term =0 for $k$ large
enough when  $t>0$ is fixed). \\ The present paper determines
$\Lambda_R$ in Section 2.
\section{Result}
Let $X$ and $Y$ be two independent random variables
following the Bernoulli $\mathcal{B}(r)$ distribution with
expectation $r$ and the Gamma $\gamma(a,\lambda)$ distribution
respectively
$$P(X=0)=1-r \hbox{ , }
 P(X=1)=r \in (0,1),$$ and $$\gamma(a,\lambda)(dt)=\frac{{\lambda}^a}{\Gamma(a)}
 t^{a-1} e^{-\lambda t}\textbf{1}_{(0,+\infty)}(t)dt.$$
The law of $X+Y$ can be seen as the mixture of Gamma distributions
up to translation that is
\begin{eqnarray*}
  \mu & = & \mathcal{B}(r)\ast \gamma(a,\lambda) \\
      & = & (1-r)\gamma(a,\lambda)+ r \gamma(a,\lambda)\ast \delta_1,\\
\end{eqnarray*}
\noindent \vskip-0.3cm  where $\delta_1$ denotes the Dirac measure at 1.
\\Now, we state our main result. The following
statement determines the set of $x>0$ such that $\mu^{{\ast}x}$
exists, or equivalently, in finding the set $\Lambda_R$ of
$(x,y)'s$ such that
$$t \longmapsto
\frac{\lambda^{y}(1-r)^{x}}{\Gamma(y)} \dis\sum_{k \geq
0}\frac{x(x-1)(x-2)....(x-k+1)}{k!} R^{k} e^{-\lambda
t}(t-k)_{+}^{y-1}=f(t)$$ is a positive function for all $t$  .
{\begin{thm}\label{10}  a) If $R\leq 1,$ then
$\Lambda_{R}=(0,+\infty)\times [1,+\infty).$
\\ b) If $R>1,$ then  $\Lambda_R=\mathbb{N}\times (0,+\infty).$
\end{thm}}
\begin{pf}
For $R \leq 1:$ for fixed $x>0$, we define the positive integer
$k_0=k_0(x)$ by $k_0-1\leq x < k_0$ and $$f_n(t)=
\frac{\lambda^{y}(1-r)^{x}}{\Gamma(y)} \dis\sum_{k=
0}^{n}\frac{x(x-1)(x-2)....(x-k+1)}{k!} R^{k} e^{-\lambda
t}(t-k)_{+}^{y-1}$$ and note that $f_n(t)>0$ for all $t>0$ and for
all $n\leq k_0.$ For seeing $\Lambda_R \subset (0,+\infty)\times
[1,+\infty)$ assume $(x,y)\in \Lambda_{R}$ but $y<1$ and consider
$f$ in the interval $(k_0+1,k_0+2)$
$$f(t)=\frac{\lambda^{y}(1-r)^{x}}{\Gamma(y)}(x-k_0 )\frac{x(x-1)(x-2)....(x-k_0+1)}{(k_0+1)!} R^{k_0+1}e^{-\lambda t} (t-k_0-1)^{y-1}+f_{k_0}(t).$$
Since $f_{k_0}(t)>0$ and $(x-k_0)<0$ we have the contradiction
$\dis\lim_{t\rightarrow k_0+1}f(t)=-\infty.$ To show that
$\Lambda_R \supset (0,+\infty)\times [1,+\infty)$, we fix $x>0$ and
$y\geq1$ and we show that $f(t)\geq0$ for all $t.$ This is already
true for $t< k_0+1$ since $f(t)=f_{k_0}(t)$ in that case. For
$k_0+1\leq k_1\leq t < k_1+1$ where $k_1$ is an integer, we use
the alternate series trick: consider the positive finite sequence
$(u_k)_{k=k_0}^{k_1} $ defined by
$$u_{k}=(-1)^{k-k_0} \frac{x(x-1)(x-2)....(x-k+1)}{k!} R^{k}
(t-k)^{y-1}$$ which is decreasing since for $k_0\leq k < k_1$ we
have
$$\frac{u_{k+1}}{u_{k}}=R \times \frac{k-x}{k+1}\Big{(}\frac{t-k-1}{t-k}\Big{)}^{y-1}<1$$
and thus
$$f(t)=f_{k_0-1}(t)+ \frac{\lambda^{y}(1-r)^{x}}{\Gamma(y)}\sum_{k=k_0}^{k_1} (-1)^{k-k_0} u_{k}\geq
0.$$ ii) For $R>1:$ of course $\Lambda_R \supset \mathbb{N}\times
(0,+\infty)$ is trivial. To prove that $\Lambda_R \subset
\mathbb{N}\times (0,+\infty)$, suppose that there exists $(x,y)\in
\Lambda_R $ such that $x$ is not an integer. Since \\$z
\longmapsto \Big{(}1+\frac{r}{1-r} e^z\Big{)}^x
\Big{(}\dis\frac{\lambda}{\lambda-z}\Big{)}^{y}$ is real analytic
on $(-\infty,\lambda)$ it is analytic on the strip
$S=(-\infty,\lambda)+i \mathbb{R}.$ However, since $R>1,$ or
$\log(\frac{1-r}{r})< \lambda$ this implies that $z \mapsto
1+\frac{r}{1-r} e^{z}$ has a zero  in the strip $S$ namely $z_0=i
\pi + \log(\frac{1-r}{r}).$ But the fact that $x$ is not an
integer prevents $z \mapsto \Big{(} 1+\frac{r}{1-r} e^{z}\Big{)}^x
$ from being analytic on $z_0$ and we get the desired
contradiction.
\end{pf}
\textbf{Acknowledgments} \vskip0.2cm We sincerely thank the Editor
and a referee for valuable suggestions and comments.
 \vskip0.2cm
\bibliographystyle{amsplain}

\begin{thebibliography}{00}
\bibitem{B} Barndorff-Nielsen, O.:  \textit{Information and Exponential
Families in Statistical Theory}, Wiley, New York (1978)\\
\bibitem {feller} Feller, W.:  \textit{An Introduction to Probability
Theory and Its Applications}, Wiley, New York (1971) \\
\bibitem{JOR} J{\o}rgensen, B.: \textit{The Theory of Dispersion Models}, London, Chapman \& Hall (1997)\\
\bibitem {R} Letac, G., Malouche, D., Maurer, S.:  The real powers of the convolution of a negative
Binomial distribution and a Bernoulli distribution.
\textit{Proc. A.M.S}.
\textbf{130}, 2107-2114 (2002)\\
\end{thebibliography}

\end{document}